\documentclass[12pt,twoside]{article}
\usepackage{amssymb}
\font\teneufm=eufm10 scaled \magstep1
\font\seveneufm=eufm7 scaled \magstep1
\font\fiveeufm=eufm5 scaled \magstep1
\newfam\eufmfam
\textfont\eufmfam=\teneufm
\scriptfont\eufmfam=\seveneufm
\scriptscriptfont\eufmfam=\fiveeufm

\newfam\msbfam
\font\tenmsb=msbm10 scaled \magstep1  \textfont\msbfam=\tenmsb
\font\sevenmsb=msbm7 scaled \magstep1 \scriptfont\msbfam=\sevenmsb
\font\fivemsb=msbm5 scaled \magstep1  \scriptscriptfont\msbfam=\fivemsb

\def\CC{{\mathbb C}}

\def\NN{{\mathbb N}}
\def\ZZ{{\mathbb Z}}

\def\PP{{\mathbb P}}

 \def\HollowBoxx #1#2#3{{\dimen0=#1 \advance\dimen0 by -#2
       \dimen1=#1 \advance\dimen1 by #3
        \vrule height 0pt depth #3 width #2
       \hskip -#3
       \vrule height #1 depth #3 width #3}}
 \def\LeftContraction{\mathord{\kern1.45pt \HollowBoxx{6pt}{3.5pt}{.4pt}}\,}

 \def\HollowBox #1#2#3{{\dimen0=#1 \advance\dimen0 by -#3
       \dimen1=#1 \advance\dimen1 by #3
        \vrule height #1 depth #3 width #3
        \vrule height 0pt depth #3 width #2
        \hskip -#3}}
 \def\RightContraction{\mathord{\, \HollowBox{6pt}{3.1pt}{.4pt}} \kern1.6pt}

\def\qed{{\hfill $\Box$}}
\newtheorem{theorem}{THEOREM}[section]

\newtheorem{remark}[theorem]{Remark}

\begin{document}

\begin{center}
{\Large \bf Characterization of the Unit Ball in ${\bf C}^n$ 
\medskip\\
Among Complex Manifolds of Dimension $n$}
\medskip\\
\normalsize A. V. Isaev

\end{center}

\begin{quotation} \small \sl We show that if the group of holomorphic automorphisms of a connected complex manifold $M$ of dimension $n$ is isomorphic as a topological group equipped with the compact-open topology to the automorphism group of the unit ball $B^n\subset\CC^n$, then $M$ is biholomorphically equivalent to either $B^n$ or $\CC\PP^n\setminus\overline{B^n}$.  
\end{quotation}

\pagestyle{myheadings}
\markboth{A. V. Isaev}{Characterization of the Unit Ball}

\setcounter{section}{0}
\section{Introduction}

For a complex manifold $M$ denote by $\hbox{Aut}(M)$ the group of holomorphic automorphisms of $M$. Equipped with the compact-open topology, $\hbox{Aut}(M)$ is a topological group. We are interested in characterizing complex manifolds by their automorphism groups.

One manifold that has been enjoying much attention in this respect is the unit ball $B^n\subset\CC^n$ for $n\ge 2$. Starting with the famous theorems of Wong \cite{W} and Rosay \cite{R} many results characterizing $B^n$ in terms of its automorphism group have been obtained. We mention proofs of Rosay's theorem by means of invariant metrics \cite{Kl}, by means of scaling \cite{P}, by means of analyzing the structure of the ring of holomorphic functions \cite{KK2}, as well as extensions of the theorem to the case of unbounded domains \cite{E}, domains in complex manifolds \cite{GKK} and domains (both bounded and unbounded) in infinite-dimensional complex space \cite{KK1}, \cite{BGK}, \cite{KM}. Rosay's theorem implies, in particular, that a bounded homogeneous domain in $\CC^n$ with $C^2$-smooth boundary is biholomorphically equivalent to $B^n$. A characterization result similar in spirit, but utilizing only the isotropy subgroup of a point in a complex manifold was obtained in \cite{GK}. More information on results of this kind can be found in the survey \cite{IKra1}.

Among Kobayashi-hyperbolic manifolds, $B^n$ can also be characterized as the manifold whose automorphism group has the largest dimension. Namely, if a connected complex manifold $M$ of dimension $n$ is hyperbolic, $\hbox{Aut}(M)$ admits the structure of a (real) Lie group of dimension not exceeding $n^2+2n$, and $\hbox{dim}\,\hbox{Aut}(M)=n^2+2n$ if and only if $M$ is biholomorphically equivalent to $B^n$ \cite{Ka}, \cite{Ko} (for a generalization of the latter result see \cite{IKra2}). These facts imply the following characterization of $B^n$ in the class of all hyperbolic manifolds: if $M$ is a connected complex hyperbolic manifold of dimension $n$ and the groups $\hbox{Aut}(M)$ and $\hbox{Aut}\left(B^n\right)$ are isomorphic as topological groups, then $M$ is biholomorphically equivalent to $B^n$. Indeed, since both $\hbox{Aut}(M)$ and $\hbox{Aut}\left(B^n\right)$ are Lie groups in the compact-open topology, the topological group isomorphism between these groups is in fact a Lie group isomorphism, and thus $\hbox{dim}\,\hbox{Aut}(M)=n^2+2n$.    

In the present paper we obtain an almost identical result without the assumption that $M$ is hyperbolic thus proving that $\hbox{Aut}\left(B^n\right)$ amost completely characterizes $B^n$ among {\it all}\, connected complex manifolds of dimension $n$.

\begin{theorem}\label{main}\sl Let $M$ be a connected complex manifold of dimension $n$. Assume that $\hbox{Aut}(M)$ and $\hbox{Aut}\left(B^n\right)$ are isomorphic as topological groups equipped with the compact-open topology. Then $M$ is biholomorphically equivalent to either $B^n$ or $\CC\PP^n\setminus\overline{B^n}$.
\end{theorem}

Our proof of Theorem \ref{main} in Section \ref{proof} is based on the classification of all connected $n$-dimensional complex manifolds admitting effective actions of the unitary group $U_n$ by biholomorphic transformations obtained in \cite{IKru} for $n\ge 2$ (see (\ref{classification})). Indeed, since $\hbox{Aut}(M)$ and $\hbox{Aut}\left(B^n\right)$ are isomorphic and $\hbox{Aut}\left(B^n\right)$ contains $U_n$, the group $U_n$ acts effectively on $M$, and therefore, for $n\ge 2$, $M$ is biholomorphically equivalent to one of the manifolds listed in (\ref{classification}). Hence in order to prove Theorem \ref{main} it remains to show that the only manifolds on the list whose automorphism groups are isomorphic to $\hbox{Aut}\left(B^n\right)$ are $B^n$ and $\CC\PP^n\setminus\overline{B^n}$.

\section{Proof of Theorem \ref{main}}\label{proof}

The theorem clearly holds for $n=1$, and from now on we assume that $n\ge 2$. As we noted in the introduction, since the groups $\hbox{Aut}(M)$ and $\hbox{Aut}\left(B^n\right)$ are isomorphic, the manifold $M$ admits an effective action of the unitary group $U_n$ and hence, by the classification obtained in \cite{IKru}, is biholomorphically equivalent to one of the following manifolds
\begin{equation}
\begin{array}{lll}
\hbox{(i)} & B^n, &\\
\hbox{(ii)} & \CC^n, &\\
\hbox{(iii)} & \CC\PP^n, &\\
\hbox{(iv)} & \widehat{B^n}/\ZZ_m, & \\
\hbox{(v)} & \widehat{\CC^n}/\ZZ_m, & \\ 
\hbox{(vi)} & \widehat{\CC\PP^n}/\ZZ_m, & \\
 \hbox{(vii)} & \left(B^n\setminus\{0\}\right)/\ZZ_m, &\\ 
\hbox{(viii)} & \left(\CC^{n}\setminus\{0\}\right)/\ZZ_m, & \\
\hbox{(ix)} & \left(\CC\PP^n\setminus\{0\}\right)/\ZZ_m, &\\
\hbox{(x)} & \left(B_r^n\setminus \overline{B^n}\right)/\ZZ_m, & \hbox{for some $r>1$},\\
\hbox{(xi)} & \left(\CC^n\setminus \overline{B^n}\right)/\ZZ_m, &\\
\hbox{(xii)} & \left(\CC\PP^n\setminus\overline{B^n}\right)/\ZZ_m, & \\
\hbox{(xiii)} & M_d^n/\ZZ_m, & \hbox{for some $d\in\CC^*$, $|d|\ne1$},
\end{array}\label{classification}
\end{equation}
where $m\in\NN$ and in (iv)--(xii) we have $m=|nk+1|$ for some $k\in\ZZ$, in (xiii) we have $(n,m)=1$, $\widehat{B^n}$, $\widehat{\CC^n}$ and $\widehat{\CC\PP^n}$ denote the blow-ups at the origin of $B^n$, $\CC^n$ and $\CC\PP^n$, respectively, $B^n_r$ denotes the ball of radius $r$ in $\CC^n$, $M_d^n$ denotes the Hopf manifold obtained by identifying every point $z\in\CC^{n}\setminus\{0\}$ with $d\cdot z$, and $\CC^n$ is realized in $\CC\PP^n$ as the collection of points with homogeneous coordinates $(1:z_1:\dots:z_n)$ (in particular, the origin is the point with homogeneous coordinates $(1:\overbrace{0:\dots:0}^{\hbox{\tiny $n$ times}})$).

Since $\hbox{Aut}(M)$ and $\hbox{Aut}\left(B^n\right)$ are isomorphic as topological groups, $\hbox{Aut}(M)$ has the structure of a Lie group of dimension $n^2+2n$ isomorphic to $\hbox{Aut}\left(B^n\right)$. In particular, $\hbox{Aut}(M)$ does not admit the structure of a complex Lie group. Indeed, suppose that there is a complex Lie group structure on $\hbox{Aut}\left(B^n\right)$. Since $\hbox{Aut}\left(B^n\right)$ is isomorphic to $G:=SU(n,1)/Z$, where $Z$ is the center of $SU(n,1)$, it follows that the group $G$ admits the structure of a complex Lie group. Denote by $G^c$ the group $G$ equipped with this complex structure. Since $G$ is simple, so is $G^c$ and therefore $G^c$ is the complexification of its maximal compact subgroup. However, any maximal compact subgroup of $G^c$ is isomorphic to $U_n$ and hence has dimension $n^2$. Therefore, maximal complex subgroups in $G^c$ have dimension $(n^2+2n)/2$ only for $n=2$, in which case we have $\hbox{dim}_{\CC}G^c=4$. Since the Lie algebra of $G^c$ is a complex simple Lie algebra, it is isomorphic to one of the following algebras: ${\mathfrak{sl}}_m$, $m\ge 2$, ${\mathfrak o}_m$, $m\ge 7$, ${\mathfrak{sp}}_{2m}$, $m\ge 2$, ${\mathfrak e}_6$, ${\mathfrak e}_7$, ${\mathfrak e}_8$, ${\mathfrak f}_4$, ${\mathfrak g}_2$. It is easy to check, however, that none of these algebras has dimension 4. Thus, we have shown that $\hbox{Aut}(M)$ does not admit the structure of a complex Lie group.

On the other hand, (iii) and manifolds of the forms (vi) and (xiii) are compact and hence their automorphism groups are complex Lie groups. Thus, $M$ is not holomorphically equivalent to any of these manifolds. 

Further, the automorphism groups of (ii) and manifolds of the form (viii) are not locally compact and thus cannot be isomorphic to a Lie group. 

Next, the automorphism group of a manifold of either of the forms (vii), (x), (xi) is the group $U_n/\ZZ_m$ (isomorphic to $U_n$ by means of the map $A\ZZ_m\mapsto (\hbox{det}A)^k\cdot A$). The same holds for a manifold of the form (xii) with $m>1$. However, the automorphism group of manifold (xii) for $m=1$ is $PSU(n,1)$ and thus is isomorphic to $\hbox{Aut}\left(B^n\right)$.

Further, a manifold of either of the forms (iv), (v) contains a unique copy of $\CC\PP^{n-1}$ (the one used to blow-up either $B^n$ or $\CC^n$ at the origin). Hence, every automorphism of such a manifold preserves it. This implies that the automorphism group of manifold (iv) is again $U_n/\ZZ_m\buildrel{\sim}\over{=} U_n$ and that of manifold (v) is $GL_n(\CC)/\ZZ_m$. Neither of these groups is isomorphic to $\hbox{Aut}\left(B^n\right)$.

Next, consider a manifold of the form (ix). For $m=1$ its automorphism group is the group of all matrices of the form
$$
\left(
\begin{array}{ll}
1 & a\\
0 & A
\end{array}
\right),
$$
with $a\in\CC^n$, $A\in GL_n(\CC)$, and $\ZZ_m$ is embedded in $GL_n(\CC)$ in the standard way as a subgroup of scalar matrices. For an arbitrary $m\ge 1$ the automorphism group of manifold (ix) is isomorphic to $GL_n(\CC)\ltimes H_{n,m}$, where $H_{n,m}$ is the group of homogeneous polynomials of degree $m$ of $n$ variables.\footnote{We are grateful to A. Huckleberry for this remark.}  

Thus, we have shown that $\hbox{Aut}(M)$ is not isomorphic to $\hbox{Aut}\left(B^n\right)$ unless $M$ is biholomorphically equivalent to either $B^n$ or $\CC\PP^n\setminus\overline{B^n}$. The proof is complete.

\qed    

\begin{remark}\label{all}\rm In \cite{IKru} we obtained a characterization of $\CC^n$ by an argument that was also based on classification (\ref{classification}). One should note, however, that not every manifold on list (\ref{classification}) is characterized by its automorphism group in the same way as $B^n$, $\CC\PP^n\setminus\overline{B^n}$ and $\CC^n$ are. For example, for every $r>1$ we have $\hbox{Aut}\left(B_r^n\setminus \overline{B^n}\right)=U_n$, but $B_{r_1}^n\setminus \overline{B^n}$ and $B_{r_2}^n\setminus \overline{B^n}$ are not holomorphically equivalent if $r_1\ne r_2$. The latter follows, for example, from the fact that two hyperbolic Reinhardt domains are biholomorphically equivalent if and only if they are equivalent by means of an elementary algebraic map \cite{Kru}, \cite{Sh}.  
\end{remark}

{\obeylines
Department of Mathematics
The Australian National University
Canberra, ACT 0200
AUSTRALIA
E-mail: alexander.isaev@maths.anu.edu.au
}


\begin{thebibliography}{ABCDE}

\bibitem[BGK]{BGK}ÊByun, J., Gaussier, H. and Kim, K.-T., Weak-type normal families of holomorphic mappings in Banach spaces and characterization of the Hilbert ball by its automorphism group, {\it J. Geom. Anal.} 12 (2002), 581--599.

\bibitem[E]{E} Efimov, A., Extension of the Wong-Rosay theorem to the unbounded case (translated from Russian), {\it Sb. Math.} 186 (1995), 967--976.

\bibitem[GKK]{GKK} Gaussier, H., Kim, K.-T. and Krantz, S. G., A note on the Wong-Rosay theorem in complex manifolds, {\it Complex Var. Theory Appl.} 47 (2002), 761--768.

\bibitem[GK]{GK} Greene, R. E. and Krantz, S. G., Characterization of complex manifolds by the isotropy subgroups of their automorphism groups, {\it Indiana Univ. Math. J.} 34 (1985), 865--879.

\bibitem[IKra1]{IKra1} Isaev, A. V. and Krantz, S. G., Domains with non-compact automorphism group: a survey, {\it Advances in Math.} 146 (1999), 1--38.

\bibitem[IKra2]{IKra2} Isaev, A. V. and Krantz, S. G., On the automorphism groups of hyperbolic manifolds, {\it J. Reine Angew. Math.} 534 (2001), 187--194.

\bibitem[IKru]{IKru} Isaev, A. V. and Kruzhilin, N. G., Effective actions of the unitary group on complex manifolds, {\it Canad. J. Math.} 54 (2002), 1254--1279.

\bibitem[Ka]{Ka} Kaup, W., Reelle Transformationsgruppen und invariante Metriken auf komplexen R\"aumen, {\it Invent. Math.} 3 (1967), 43--70. 

\bibitem[KK1]{KK1} ÊKim, K.-T. and Krantz, S. G., Characterization of the Hilbert ball by its automorphism group, {\it Trans. Amer. Math. Soc.} 354 (2002), 2797--2818.

\bibitem[KK2]{KK2} Kim, K.-T. and Krantz, S. G., Some new results on domains in complex space with non-compact automorphism group, {\it J. Math. Anal. Appl.} 281 (2003), 417--424.

\bibitem[KM]{KM} ÊKim, K.-T. and Ma, D., Characterization of the Hilbert ball by its automorphisms, {\it J. Korean Math. Soc.} 40 (2003), 503--516.

\bibitem[Kl]{Kl} Klembeck, P., K\"ahler metric of negative curvature, the Bergman metric near the boundary, and the Kobayashi metric on smooth bounded strictly pseudoconvex sets, {\it Indiana Univ. Math. J.} 27 (1978), 275--282.

\bibitem[Ko]{Ko} Kobayashi, S., {\it Hyperbolic Manifolds and Holomorphic Mappings}, Marcel Dekker. New York 1970.

\bibitem[Kru]{Kru} Kruzhilin, N. G., Holomorphic automorphisms of hyperbolic Reinhardt domains (translated from Russian), {\it Math.\ USSR-Izv.} 32(1989),
15--38.

\bibitem[P]{P} Pinchuk, S. I., Holomorphic mappings in $\CC^n$ and the problem of holomorphic equivalence (translated from Russian),
{\it Several complex variables. III. Geometric function theory, Encycl. Math. Sci.} 9 (1989), 173-199.

\bibitem[R]{R} Rosay, J. P., Sur une charact\'erisation de la boule parmi les domaines de $\CC^n$ par son groupe d'automorphismes, {\it Ann. Inst. Fourier (Grenoble)} 29 (1979), 91--97.

\bibitem[Sh]{Sh} Shimizu, S., Automorphisms of bounded Reinhardt domains, {\it Japan. J. Math.} 15(1989), 385--414.

\bibitem[W]{W} Wong, B., Characterization of the unit ball in $\CC^n$ by its automorphism group, {\it Invent. Math.} 41 (1977), 253--257.






\end{thebibliography}
\end{document}